# Stability of Oscillator Ising Machines:
# Not All Solutions Are Created Equal

Mohammad Khairul Bashar, Zongli Lin, Nikhil Shukla*

Department of Electrical and Computer Engineering, University of Virginia,

Charlottesville, VA- 22904 USA

*e-mail: ns6pf@virginia.edu

**Abstract:**

Nonlinear dynamical systems such as coupled oscillators are being actively investigated as Ising machines for solving computationally hard problems in combinatorial optimization. Prior works have established the equivalence between the global minima of the Lyapunov function (commonly referred to as the energy function) describing the coupled oscillator system and the ground state of the Ising Hamiltonian. However, the properties of the oscillator Ising machine (OIM) from a nonlinear control viewpoint, such as the stability of the OIM solutions remains unexplored. Therefore, in this work, using nonlinear control-theoretic analysis, we (i) Identify the conditions required to ensure the functionality of the coupled oscillators as an Ising machine; (ii) Show that all globally optimal phase configurations may not always be stable, resulting in some configurations being more favored over others, and thus, creating a biased OIM; (c) Elucidate the impact of the stability of locally optimal phase configurations on the quality of the solution computed by the system. Our work, fostered through the unique convergence between nonlinear control theory and analog systems for computing, provides a new toolbox for the design and implementation of dynamical system-based computing platforms.

Oscillator Ising machines (OIMs) provide an elegant dynamical system platform for minimizing the Ising Hamiltonian $H = -\sum W_{ij} s_i s_j$, where $s_i \in \{1, -1\}$ corresponds to the $i^{th}$ spin, and $W_{ij}$ is the interaction coefficient between nodes $i$ and $j$; the Zeeman term has been neglected here. Minimizing $H$ is a quintessential NP-hard combinatorial optimization problem (COP). Moreover, such systems are being actively investigated for solving many computationally challenging COPs which can be directly mapped to the minimization of $H$[1,2]. An archetypal example of such a mapping is the MaxCut problem which is defined as a graph cut that maximizes the number of cut edges (unweighted graph considered here). The relationship between $H$ and the MaxCut ($MC$) can be defined as $H = \Sigma W - 2MC$, where $\Sigma W$ is the sum of the weights of all the edges in the graph. Thus, the optimal MaxCut directly corresponds to the ground state (minimum $H$) of a topologically equivalent spin network with antiferromagnetic interactions i.e., $W_{ij} = -1$. The promise of OIMs is that allowing the physics to do the computation can potentially provide a significant benefit in computational performance over digital algorithms[3,4].

In groundbreaking work by Wang *et al.*[5], the authors demonstrated that a global minimum of the Lyapunov function for a topographically equivalent coupled oscillator network under second harmonic injection is equivalent to computing a global minimum of $H$. While the minimization of the energy function in OIMs[5] as well as their implementation[3,6–13] has been explored in prior work, the stability of the globally optimal and locally optimal minima (attractors) and the resulting impact on the OIM dynamics has been largely unexplored. The works by Erementchouk *et al.*[14] and Böhm *et al.*[15] are the few examples that aim to investigate the dynamics of the OIM while a few more works have focused on analyzing the dynamical properties in spiking neural network[16–19]. Consequently, understanding the properties of the OIM as a nonlinear dynamical system and elucidating their impact of the computational properties are the primary focus of this work.

The dynamics of the OIM are such that the oscillator phases settle to $\theta \in \{0, \pi\}$, which subsequently, represent $s = \pm 1$ assignment to the nodes. The computational capability of this system arises from the fact that the resulting phase configuration of the oscillators will correspond to a ground state of $H$. The Lyapunov function $E(\theta(t))$ and the corresponding system dynamics are respectively presented as

$$E(\theta(t)) = -K \sum_{i,j=1, j \neq i}^{N} W_{ij} \cos(\theta_i - \theta_j) - K_s \sum_{i=1}^{N} \cos(2\theta_i(t)) \tag{1}$$

and

$$\frac{d\theta_i(t)}{dt} = -K \sum_{j=1, j \neq i}^{N} W_{ij} \sin(\Delta\theta_{ij}) - K_s \sin(2\theta_i(t)) \tag{2}$$

where [W] represents the coupling matrix between nodes, and $K$ and $K_s$ represent the strength of coupling among the oscillators and the strength of the second harmonic injection signal, respectively. Using equations (1) and (2), it can be shown that $\frac{dE(\theta)}{dt} = -2 \sum_{i=1}^{N} \left(\frac{d\theta_i(t)}{dt}\right)^2 \leq 0$,[5,20] which consequently, implies that the system will evolve towards the ground state except when $\frac{dE(\theta)}{dt} = \frac{d\theta(t)}{dt} = 0$. A point in the phase space where $\frac{d\theta(t)}{dt} = (-\nabla E) = 0$ defines a local energy minimum and there are multiple such points in the phase space. In fact, every possible spin assignment and its equivalent in terms of the oscillator phases $\{\theta_1, \theta_2, \dots, \theta_N\}$, where $\theta_i \in \{0, \pi\}$, can correspond to a local energy minimum. Consequently, the phase space contains $2^N$ minima in the system; $2^{N-1}$ minima when symmetricity in the solutions is considered. We refer to the energy minima lying at the lowest possible energy as globally optimal minima while the others are referred to as the locally optimal minima. We show that even though a point in the phase space may be a local or global minima where $\frac{d\phi(t)}{dt} = 0$, the system dynamics may not always settle or be stable at those points. Furthermore, the system dynamics can preferably converge and stabilize in some energy minima (including globally optimal solutions) over others lying at the same energy. This implies that the system may intrinsically favor certain solutions over others leading to a *biased OIM*. Consequently, engineering the system stability can have significant impact on the computational characteristics and the performance of the system.

To elucidate our approach, we consider an illustrative randomly generated unweighted graph with 20 nodes and 152 edges as shown in Fig. 1a. Figure 1b shows a histogram for the energy (quantified using $H$ here) for all possible solutions. It can be observed in Fig. 1b that the graph has 22 spin configurations that yield the minimum energy ($H = -28$). However, as alluded to above, the system dynamics may not always be stable for all the 22 globally optimal configurations.

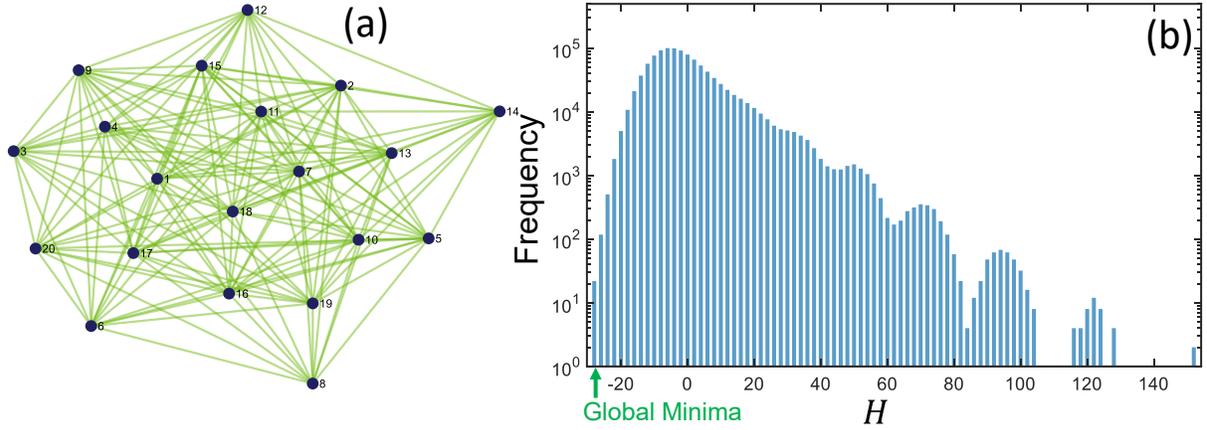

**Fig. 1.** (a) Illustrative randomly graph with 20 nodes and 152 edges. (b) Corresponding histogram of energy ($H$) for all ($2^{20}$) possible spin configurations.

In order to investigate the stability of the global and local solutions, we analyze the Lyapunov exponents $(\lambda_1, \lambda_2, \lambda_3, \ldots, \lambda_N)$ (eigenvalues of the Jacobian matrix[21]) for the system dynamics. Lyapunov exponents provide a powerful mathematical tool for analyzing the stability of non-linear dynamical systems[22]. In order for a phase configuration to be stable, no Lyapunov exponent should be greater than zero. The Jacobian matrix ($J$) for the OIM (assuming symmetric unweighted edges i.e., $W_{ij} = W_{ji}$) can be defined as,

$$J = \begin{bmatrix} E(1,1) & -KJ_{12}\cos(\theta_1 - \theta_2) & -KJ_{13}\cos(\theta_1 - \theta_3) & \ldots & -KJ_{1N}\cos(\theta_1 - \theta_N) \\ -KJ_{12}\cos(\theta_1 - \theta_2) & E(2,2) & -KJ_{23}\cos(\theta_2 - \theta_3) & \ldots & \\ -KJ_{13}\cos(\theta_1 - \theta_3) & -KJ_{23}\cos(\theta_2 - \theta_3) & E(3,3) & \ldots & \vdots \\ \vdots & \vdots & \vdots & & \\ -KJ_{1N}\cos(\theta_1 - \theta_N) & & \ldots & & E(N,N) \end{bmatrix} \quad (3)$$

where $E(i,i) = -K \sum_{j=1, j \neq i}^{N} W_{ij} \cos(\theta_i - \theta_j) - 2K_s \cos(2\theta_i)$. The Eigenvalues of $J$ for a given point in the phase space where $\frac{d\phi(t)}{dt} = 0$ yield the Lyapunov exponents at that point. Since *all* the Lyapunov coefficients need to be negative in order for an energy minimum to be stable, we focus on the largest Lyapunov coefficient (referred to here as $\lambda_L$) since all other coefficients will be smaller than $\lambda_L$.

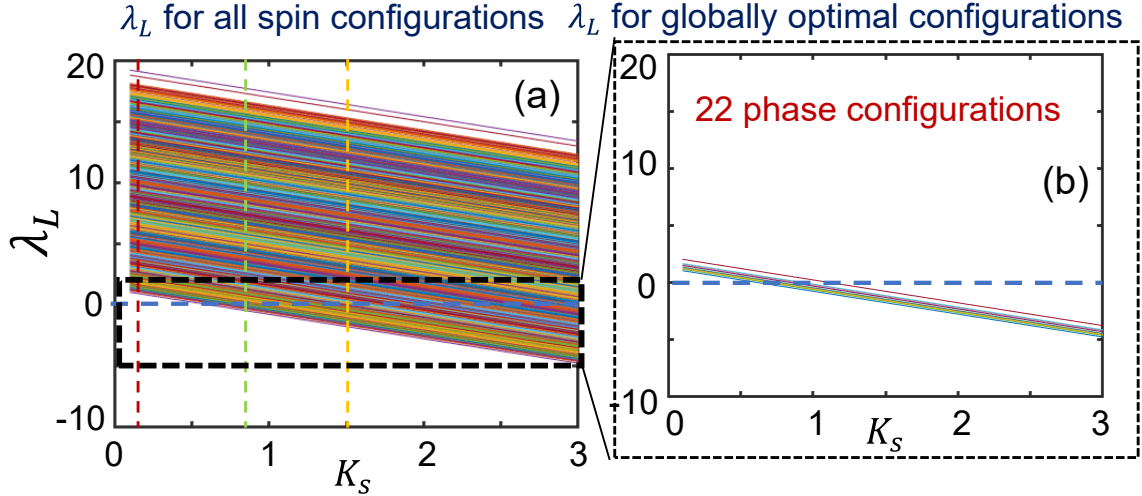

**Fig. 2.** (a) Evolution of the largest Lyapunov coefficient ($\lambda_L$) as a function of $K_s$ for all spin configurations; (b) Evolution of $\lambda_L$ as a function of the $K_s$ for the subset of globally optimal phase configurations. Note: $\lambda_L > 0$ implies that the particular solution is unstable.

Figure 2a shows the evolution of the *largest* Lyapunov exponent ($\lambda_L$) as a function of $K_s$ ($K = 1$) for the representative graph shown in Fig. 1a. All the $2^{20}$ possible phase configurations are considered. The evolution of $\lambda_L$ for only the globally optimal solutions is emphasized in Fig. 2b. It can be observed that the stability of a spin configuration is significantly impacted by the strength $K_s$ (relative to $K$) of the second harmonic injection signal. In fact, if $K_s$ is small enough (<0.5 for the graph considered here), then the ground states, i.e., globally optimal configurations themselves can become unstable. In such a scenario, the system will cease to behave as an Ising machine – the ground state energy of the system will then correspond to an oscillator phase configuration where some or all oscillator phases do not settle to 0 or $\pi$.

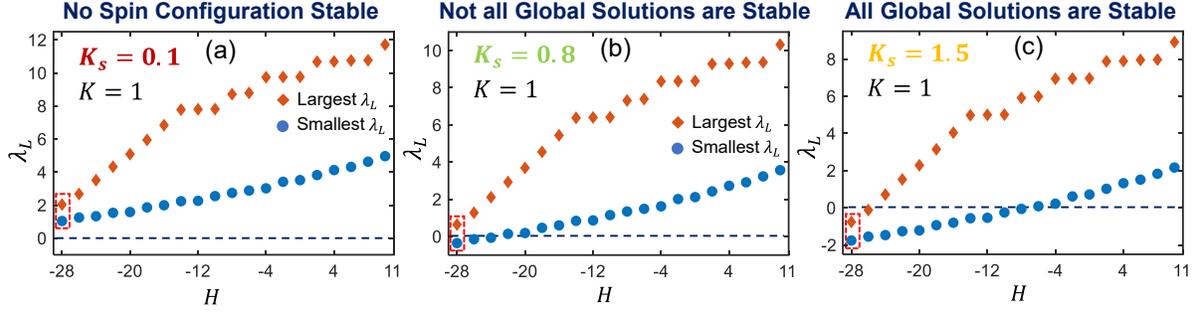

**Fig. 3.** Minimum (blue) and maximum (orange) $\lambda_L$ for phase configurations lying at a particular energy ($H$). Three values of $K_s$ are considered (a) $K_s = 0.1$: Since all spin configurations (including the globally optimal solutions) are unstable, it implies that the oscillator platform will cease to behave as an Ising machine. (b) $K_s = 0.8$: Some globally optimal solutions are stable while others are unstable. Additionally, a few locally optimal low energy solutions are also stabilized; (c) $K_s = 1.5$: All globally optimal solutions are stable. The red box indicates the globally optimal solutions.

Next, for different $K_s$ ($K = 1$), we analyze the distribution of $\lambda_L$ for all phase configurations lying at a given energy ($H$). Figs. 3a-c show the maximum and the minimum value of $\lambda_L$ for phase configurations corresponding to a given $H$, computed for three different values of $K_s$ (0.1, 0.8, 1.5), respectively. In Fig. 3a, it can be observed that since $\lambda_L$ for all spin configurations (including the globally optimal solutions lying at $H = -28$) are greater than zero, the ground state of the oscillator platform will not be achieved for $\theta \in \{0, \pi\}$. Consequently, it is expected that the oscillator platform will cease to function as an Ising machine. When $K_s = 0.8$ (Fig. 3b), it can be observed that the maximum and the minimum value of $\lambda_L$ for the globally optimal solutions straddle zero, i.e., $\lambda_L$ for some solutions is less than zero whereas it is greater than zero for others. This implies that only a fraction of the globally optimal spin configurations are stable, and consequently, the system dynamics will preferentially converge to the stable (globally optimal) solutions. This creates a biased OIM that favors the stable states over the unstable ones. Additionally, it is also noteworthy to point out that some of the locally optimal (but globally sub-optimal) solutions lying at low energies ($H = -24, -26$) are also stabilized. This indicates that the system may potentially get trapped in one of these states, leading to sub-optimal solutions. However, the value of $K_s$ is such that the solutions lying at higher energies ($H > -24$) are destabilized, preventing the system from getting trapped in those states. Finally, when the strength of the second harmonic injection is increased further to $K_s = 1.5$ (Fig. 3c), it can be observed that all the globally optimal solutions are stabilized. Additionally, increasing the strength of the second harmonic injection also increases the number and energy of the locally optimal (globally sub-optimal) solutions where the system dynamics can be stabilized. Consequently, this should increase the probability of the system getting trapped at a local minimum.

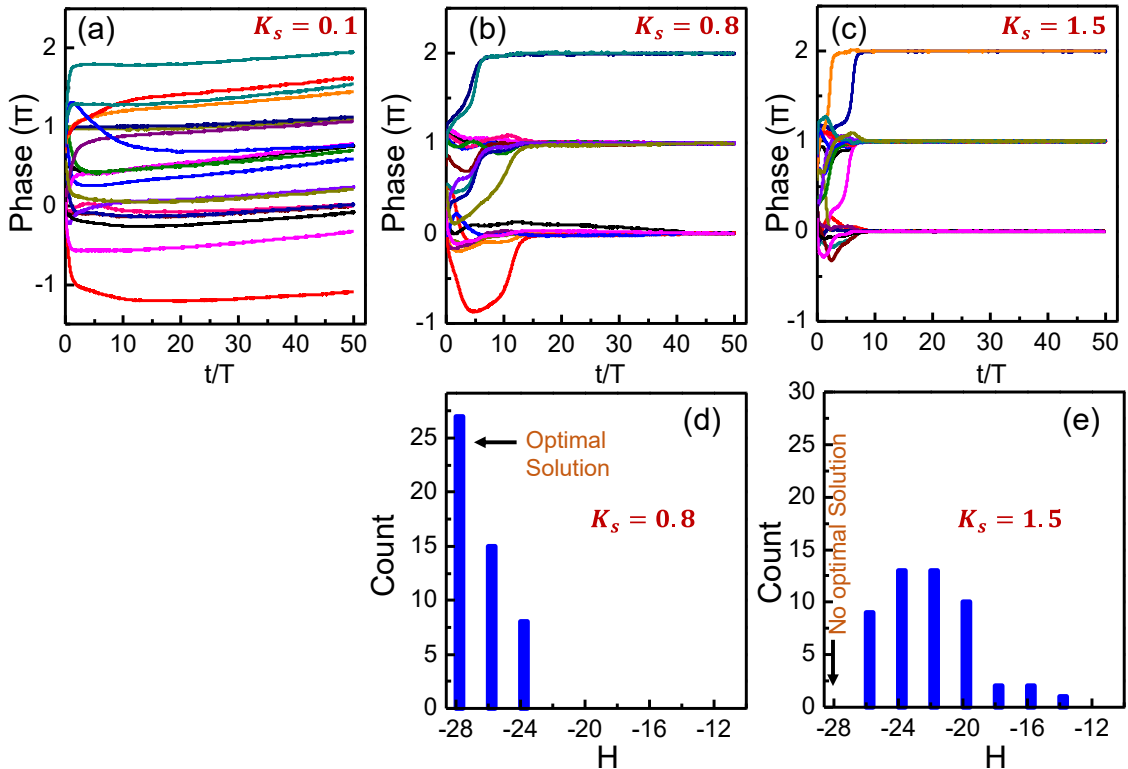

**Fig. 4.** Temporal evolution of the oscillator phases for (a) $K_s = 0.1$; (b) $K_s = 0.8$; (c) $K_s = 1.5$, respectively ($K = 1$). The oscillator network is topologically equivalent to the graph considered in Fig. 1(a). When $K_s = 0.1$, the phases do not converge to {0, π}, and thus, the system does not behave as an Ising machine. Measured $H$ for (d) $K_s = 0.8$. (e) $K_s = 1.5$ over 50 trials. Since a smaller number of locally optimal solutions are stabilized at $K_s = 0.8$ compared to $K_s = 1.5$, the system yields better solution quality. No solution is obtained for $K_s = 0.1$.

We verify the system behavior predicted above using simulations shown in Fig. 4. We consider an oscillator network that is equivalent to the graph considered in Fig. 1a, and subsequently, evaluate the dynamics for different second harmonic injection strengths. We simulate the system dynamics (2) using MATLAB's® SDE (stochastic differential equation) solver, where we consider a time and phase independent noise amplitude of $K_n = 0.005$. When $K_s = 0.1$, the oscillator phases, as expected, do not converge to $\{0, \pi\}$ and the oscillator platform does not behave as an Ising machine. For larger injection strengths ($K_s = 0.8, 1.5$), it can be observed that the oscillator phases are binarized to $\{0, \pi\}$, validating the system's ability to function as an Ising machine.

Figures 4d and 4e show a histogram of the computed $H$ for $K_s = 0.8$ and $K_s = 1.5$, respectively, over 50 trials with randomly generated initial conditions. The spin assignments and the solutions are not computed for $K_s = 0.1$ since the system does not behave as an Ising machine. Since $K_s = 0.8$ only stabilizes some globally optimal solutions and some phase configurations that lie at low energy ($H = -26, -24$), it can be observed that the system dynamics always converge to one of these states. In contrast, increasing the second harmonic injection strength to $K_s = 1.5$ stabilizes all the 22 global solutions as well as many other phase configurations that lie at higher energies (Fig. 3c). Consequently, it can be observed in Fig. 4e that the system dynamics exhibit a higher probability of getting trapped at a local minimum (sub-optimal spin configuration). In fact, over 50 trials, the system never converges to a globally optimal solution ($H = -28$). This indicates that the ability to engineer the stability of the local minima can significantly impact the computational performance of the OIM.

In summary, we have uniquely examined the OIM from a control theoretic standpoint. By calculating the Lyapunov coefficients for various spin configurations, we analyze the stability of the globally and locally optimal solutions. We show that not all globally optimal phase configurations may always be stable, in which case the system behaves as a biased OIM. Furthermore, our analysis also provides insights into the stability of the locally optimal phase configurations which can significantly impact the probability of the system dynamics getting trapped at a local minimum, which subsequently, impacts the quality of the solution computed by the system. This work presents a novel paradigm rooted in nonlinear control theory for analyzing dynamical system-based computing platforms such as OIMs and creates a new toolkit for the design and implementation of such systems.

**Acknowledgment**

This work was supported by NSF grant No. 2132918.

**Competing interests**

The authors declare no competing interests.

**Data availability**

The data that support the findings of this study are available from the corresponding author upon reasonable request.

**Author contributions**

M. K. B performed the simulations. Z. L. helped the analysis and interpretation of the results. N. S. conceived the idea, designed the research and helped with the simulations. N.S. wrote the manuscript. Z. L helped proof read the manuscript.